\documentclass[a4paper, 12pt]{article}
\usepackage{amsmath} \usepackage{amsfonts} \usepackage{amssymb}\usepackage{latexsym}
\usepackage[english]{babel}
\usepackage{amscd}
\usepackage[dvips,final]{graphics}
\usepackage{locengli,math}
\usepackage{graphicx, Addaletree}
\usepackage[all]{xy}
\begin{document}
\begin{center}
\textbf{\large{\textsf{Infinitesimal or cocommutative dipterous bialgebras \\and\\  good triples of operads}}}
\vskip1.5cm
\parbox[t]{14cm}{\large{
Philippe {\sc Leroux}}\\
{\footnotesize
ph$\_$ler$\_$math@yahoo.com}}
\end{center}

\vskip1.5cm
\noindent
{\bf Abstract:}

\noindent
\textbf{Notation}:
In the sequel $K$ is a field and $\Sigma_n$ is the group of permutation over $n$ elements. If $\mathcal{A}$ is an operad, then the $K$-vector space of $n$-ary operations is denoted as usual by $\mathcal{A}(n)$. Recall that if $\mathcal{A}$ is regular, then $\mathcal{A}(n):= \mathcal{A}_n \otimes K \Sigma_n$, where $\mathcal{A}_n$ is the $K$-vector space of $n$-ary operations without permutations of the entries. We adopt Sweedler notation for binary cooperation $\Delta$ on a $K$-vector space $V$ and set $\Delta(x)= x_{(1)} \otimes x_{(2)}$. Left dipterous algebras in the sequel will be just
abbreviated as dipterous algebras.

\section{Introduction}
The works of Poincar\'e, Birkhoff, Witt (P.B.W.) and Cartier, Milnor, Moore (C.M.N.) on the connected cocommutative Hopf algebras can be summarized as follows.
For any cocommutative bialgebra $\mathcal{H}$, the following are equivalent:
\begin{enumerate}
\item{$\mathcal{H}$ is connected;}
\item{There is an isomorphism of bialgebras $\mathcal{H} \simeq U(Prim \ \mathcal{H})$;}
\item{There is an isomorphism of connected coalgebras $\mathcal{H} \simeq Com^c(Prim \ \mathcal{H})$, }
\end{enumerate}
where $U: Lie-alg. \rightarrow As-alg.$ is the usual universal enveloping algebra functor, $Com^c(V)$ is the cofree cocommutative coassociative coalgebra over a $K$-vector space $V$ and $Prim \ \mathcal{H}$ is the Lie algebra of the primitive elements of $\mathcal{H}$.
In other words, the triple of operads $(Com, As, Lie)$ is good according to the terminology of Loday \cite{GB}. Since then
many other good triples have been found such as for instance the triple $(As,As,Vect)$ \cite{LodRon} endowed with the nonunital infinitesimal compatibility relation,
$$ \delta(xy) = x_{(1)}\otimes x_{(2)} y + xy_{(1)}\otimes y_{(2)} +x \otimes y.$$
The reader should read \cite{GB} $p.102$ for a summary.

\noindent
We focus on this paper to good triples involving the operad $Dipt$ (resp. $RDipt$) of dipterous (resp. right dipterous) algebras instead of $As$.
A dipterous algebra is a $K$-vector space equipped with two binary operations $\star$ and $\succ$
verifying:
$$ (x \star y) \star z= x \star (y \star z), \ \ (x \star y) \succ z= x \succ  (y \succ  z).$$
Similarly, a right dipterous algebra
is a $K$-vector space equipped with two binary operations $\star$ and $\prec$
verifying:
$$ (x \prec y) \prec z= x \prec (y \star z), \ \ (x \star y) \star z= x \star  (y \star  z).$$
They have been introduced by Loday and Ronco in \cite{LodRon} and also by the author in his thesis (see \cite{dipt}) via their coversions where right codipterous coalgebras\footnote{In this text, right codipterous coalgebras are named anti-codipterous coalgebras.} entangled together were the elementary bricks of our coassociative geometry over directed graphs.
We set
$Dipt-alg.$, resp. $Rdipt-alg.$, the associated categories. We mainly focus on the operad $Dipt$ since
all our results may be carried on $RDipt$ very straightforwardly.
In \cite{LodRon}, Loday and Ronco showed that the triple of operads $(As, Dipt, B_\infty)$ endowed with the semi-Hopf relations was good.
In this paper we provide other good triples involving the operad $Dipt$.
In Section 2, many examples of dipterous algebras are given, notably the free L-dipterous algebra over $V$ \cite{BaxLer} which is closely related to duplicial algebras \cite{GB}/ triplicial algebras \cite{Ltrip}.
In Section 3,
we construct explicitely the free dipterous algebra over a $K$-vector space $V$. This construction involved forests of planar rooted trees
and was announced in \cite{LodRon, LodRon1} to be so. As a corollary, we also propose a injective coding of rooted planar trees via rooted planar $m$-ary trees, $m>1$.
In Section 4, we introduced infinitesimal dipterous bialgebras as dipterous algebras equipped with a coassociative coproduct $\Delta$ verifying new compatibility relations called nonunital semi-infinitesimal relations:
$$ \Delta(x \succ y):= x_{(1)}\otimes (x_{(2)} \succ y) + (x \star y_{(1)}) \otimes y_{(2)} + x\otimes y.$$
$$ \Delta(x \star y):= x_{(1)}\otimes (x_{(2)} \star y) + (x \star y_{(1)}) \otimes y_{(2)} + x\otimes y,$$
which should not be confused with the nonunital infinitesimal relation:
$$ \Delta(x \vdash y):= x_{(1)}\otimes (x_{(2)} \vdash  y) + (x \vdash  y_{(1)}) \otimes y_{(2)} + x\otimes y,$$
which is used in \cite{GB, HLR, Ltrip} for instance, the two relations coinciding only for associative products.
We prove then that the triple of operads
$(As, Dipt, Mag^{\infty})$ is good, where the operad $Mag^{\infty}$ is explicitely described with the help of rooted planar trees in Section 5.
In Section 6, the dual in the sense of Ginzburg and Kapranov of $Dipt$, called the operad $QNDipt$ is also given. An Homology of dipterous algebras is given and the operad $Dipt$ turns out be Koszul. A 2-associative algebra is a $K$-vector space equipped with two associative products  \cite{LodRon}. In Section 7, we prove a rigidity theorem for the so-called connected $2As^c-Dipt$-bialgebras, i.e., the triple of operads $(2As, Dipt, Vect)$ endowed both with the unital semi-Hopf and with the unital semi-infinitesimal compatibility relations is good.
We close this paper by proposing another good triple of operads involving $Dipt$, related to the Connes-Kreimer Hopf algebra formalism in quantum field theory, the triple $(Com, Dipt, Prim_{Com} Dipt)$ endowed with the Hopf compatibility relations, the operad of primitive elements being
unknown and we generalize our results
to dipterous like operads and associative molecules.
\section{Examples}
By definition, any associative algebra  with right or left module over itself provide dipterous algebras. We give now other examples.
\subsection{From coassociative manifolds}
In \cite{dipt} codipterous coalgebras are constructed from coalgebras to provide constructions of directed graphs geometric supports for coassociative coproducts.

\subsection{From Language theory}
Interesting dipterous algebras can be constructed via the notion of L-dipterous algebras \cite{BaxLer}. A L-dipterous algebra $(A, \star, \succ)$ is a dipterous algebra verifying the following extra condition:
$$  \ \ \ (x \succ y) \star z = x \succ (y \star z),$$
for all $x,y \in A$.
Respectively, a (right) L-dipterous algebra $A$ is an right dipterous algebra with:
$$ \ \ (x \star y)\prec  z = x \star (y \prec z),$$
for all $x,y \in A$.

\NB
Because of the last relation, we have a functor $L-Dipt\rightarrow L$, hence the name. L-algebras have been first introduced in \cite{Coa}, see \cite{Ltrip} for more information.

\noindent
It has been shown in \cite{perorb1} that Language theory or symbolic dynamics can
viewed through the use of cooperations, that is cooperads.
The following is taken from our unpublished paper \cite{BaxLer}. In this paper are introduced  {\it{$\epsilon'[R]$-bialgebras}} (resp. {\it{$\epsilon'[L]$-bialgebras}}).
Such an object
$(A,\mu,\Delta)$ is an associative algebra $(A,\mu)$ together with a coalgebra (not necessarily coassociative) $(A,\Delta$) such that, for all $a,b \in A$, $\Delta(ab) := a \Delta(b)$, (resp. $\Delta(ab) := \Delta(a) b$).
Similarly, the notion of right Baxter-Rota operators (resp. left Baxter-Rota operators) is introduced. If $\mathcal{A}$ is a binary operad, and $A$ is a $\mathcal{A}$-algebra, then such operators are linear maps $\zeta: A \xrightarrow{} A$ verifying
$\zeta(x) \diamond \zeta(y)=\zeta(\zeta(x)\diamond y))$ (resp. $\zeta(x)\diamond \zeta(y)=\zeta(x\diamond \zeta(y))$), for all generating binary operation $\diamond$ of $\mathcal{A} $.
We now focus on $\epsilon'[R]$-bialgebras. The following results show that associative structures can pop up from noncoassociative cooperations provided they well behave with the underlying associative product.
\begin{prop}\cite{BaxLer}
\label{pop}
\begin{enumerate}
\item{Let $(A, \ \mu)$ be an associative algebra and $\zeta: A \xrightarrow{} A$ be a right Baxter-Rota operator. Define the binary operations, $\star_\zeta, \ \prec_\zeta: A^{\otimes 2} \xrightarrow{} A$ by:
$$x \star_\zeta y := \zeta(x)y, \ \ \ \ \
x \prec_\zeta y := x \zeta(y), \ \forall \ x,y \in A.$$
Then, $(A, \ \star_\zeta, \ \prec_\zeta)$ is a right $L$-dipterous algebra.}
\item{Let $(A, \mu, \Delta)$ be a $\epsilon'[R]$-bialgebra.
Equip the algebra $\textsf{End}(A)$ with the convolution product $*$.
Then, there exists a right Baxter-Rota operator, $\beta: \textsf{End}(A) \xrightarrow{} \textsf{End}(A)$ given by $\beta(T):= id * T$, for all $T \in \textsf{End}(A)$. Set,
$$T \star_\beta S := \beta(T) S \ \ \ \ \textrm{and} \ \ \ \
T \prec_\beta S := T  \beta(S), \ \forall \ T,S \in \textsf{End}(A).$$
Then, $(\textsf{End}(A), \star_\beta, \prec_\beta)$ is a right $L$-dipterous algebra.}
\item{Let $As(S)$ be the free associative $K$-algebra generated by a nonempty set $S$. Fix a cooperation $\Delta: KS \xrightarrow{} KS^{\otimes 2}$ and extend it to a cooperation $\Delta_{\sharp}: As(S) \xrightarrow{} As(S)^{\otimes 2} $ defined for any words $s_1 \ldots s_n$ by $\Delta_{\sharp}(s_1 \ldots v_s):= s_1 \ldots s_{n-1} \Delta(s_n)$.
Then, $(As(S),\Delta_{\sharp})$ is a $\epsilon'[R]$-bialgebra.}
\item{
Let $(A, \ \mu, \ \Delta)$ be a $\epsilon'[R]$-bialgebra.
Set $x \bowtie y := \mu(\Delta(x)) y$ and $x \prec_{A} y := x \mu(\Delta(y))$ defined for all $x,y \in A$.
Then, the $K$-vector space $(A, \ \bowtie, \ \prec_{A} )$ is a right $L$-dipterous algebra.
}
\end{enumerate}
\end{prop}
\begin{exam}{[Weighted directed graphs and dynamics]}
\label{graphdyn}
Let $G=(G_0, G_1, s, t)$ be a directed graph, supposed to be locally-finite, row-finite, without sink and source, equipped with a family of weights $(w_v)_{v \in G_0}$ and such that $s \times t: G_1 \xrightarrow{} G_0 \times G_0$ is injective. Consider the free $K$-vector space $KG_0$ spanned by $G_0$. Identify any directed arrow $v \longrightarrow w \in G_1$ with $v \otimes w$.
The set $G_1$ is then viewed as a subset of $KG_0^{\otimes 2}$. The family of weights $(w_v)_{v \in G_0}$ is then viewed as
a family of maps $w_v: \ F_v \xrightarrow{} k$, where $F_v:=\{a \in G_{1}, \ s(a)=v \}$.
Define the co-operation $\Delta_M: KG_0 \xrightarrow{} KG_0^{\otimes 2}$ as follows \cite{Coa}:
$$\Delta_M (v) := \sum_{i: a_i \in F_v} \ w_v(a_i) \ v \otimes t(a_i),$$ for all $v \in G_0$. Extend it to $As(G_0)$ as in Proposition \ref{pop}. Then, $(As(G_0), \Delta_{M\sharp})$ is a $\epsilon'[R]$-bialgebra.
\end{exam}
\begin{exam}{[Substitutions and Language theory]}
Let $S$ be a non-empty set. For all $s \in S$, suppose there exit substitutions of the form $s \mapsto s^i
_{1}s^i _{2}$, $i \in I_s$, which appear with probability $\mathbb{P}(s \mapsto s^i _{1}s^i _{2})$ and $\card(I_s) < \infty$ for all $s \in S$. Consider the free
$K$-vector space $K S$ and define the cooperation
$\Delta: KS  \xrightarrow{} K S^{\otimes 2}$ by,
$$\Delta (s) := \sum_{i \in I_s} \ \mathbb{P}(s \mapsto s^i _{1}s^i _{2}) \ s^i _{1} \otimes s^i _{2}.$$
Extend it to $As(S)$ as in Proposition \ref{pop}. The $K$-vector space $(As(S), \ \Delta_{\sharp})$ is then a $\epsilon'[R]$-bialgebra.
All the possible dynamics of a string are then given by the operator $\mu \Delta$, where $\mu$ is the associative product representing the concatenation of two symbols.
For instance, consider the start symbol at time $t=0$, say $s$. Apply
the operator $\mu\Delta$ to obtain,
$$\sum_{i \in I_s} \ \mathbb{P}(s \mapsto s^i _{1}s^i _{2}) \ s^i _{1} s^i _{2},$$
which is all the possible strings at time $t=1$. The probability to get the word, say $s^{i_0} _{1} s^{i_0} _{2}$, with $i_0 \in I_s$ at time $t=1$, is then
$\mathbb{P}(s \mapsto s^{i_0} _{1}s^{i_0} _{2})$ and so forth.
\end{exam}
\subsubsection{The free L-dipterous algebra}
The aim of this subsection is to construct the free L-dipterous algebra over a $K$-vector space $V$.
We consider the rooted planar binary trees and denote by $Y_n$ the set of such trees with $n$ internal degrees. For instance,
$Y_0:=\{\vert \}, Y_1:=\{ \treeA \}$,  $Y_2:=\{ \treeAB, \ \treeBA \}$.
Recall that any such trees can be uniquely written via the grafting operation $\vee$. If $r \in Y_p$ and $s \in Y_q$ then $r \vee s \in Y_{p+q+1}$ is the tree obtained from $r$ and $s$ by gluing their roots together and adding a new root: $^r \treeA ^s$. For instance $\treeBA= \vert \vee \treeA$.
A tree $t$ will be decomposed as  $t= t_l \vee t_r$.
Define the binary  operation $\nwarrow$ over trees as follows. If $r \in Y_p$ and $s \in Y_q$ then $r \nwarrow s \in Y_{p+q}$ is the tree obtained from $r$ and $s$ by gluing the root of $s$ and the most right leaf of $r$ together. For instance $\treeBA=\treeA \nwarrow \treeA$. Define also the binary  operation $\succ$ over trees as follows. If $t \in Y_p$ and $s \in Y_q$ then $t \succ s \in Y_{p+q}$ is defined by $ t \succ s := (t \nwarrow s_l)\vee s_r$. We set $\vert \nwarrow t= t\nearrow \vert=t$.
\begin{theo}
The free L-dipterous algebra over $K$ is:
$$ Y_{\infty}:= \bigoplus_{n>0} \ KY_n,$$
as a $K$-vector space equipped with the two binary operations $\nwarrow$ and $\succ$ extended by bilinearity.
\end{theo}
\Proof
The two binary operations $\nwarrow$ and $\succ$ defined just above are extended by bilinearity. The operation $\nwarrow$ is associative. Moreover, let $r,s,t$ be trees. We get $(t\nwarrow r) \succ s =(t_l\vee (t_r \nwarrow r \nwarrow s_r))\vee s_l= t\succ (r \succ s)$. Similarly, we get
$t \succ (r \nwarrow s) = (t\succ r)\nwarrow s$, proving that $Y_{\infty}$ is a L-dipterous algebra. Let $LD$ be a L-dipterous algebra and $x \in LD$. Let $f: K \rightarrow LD$, $1_K \mapsto x$ be a linear map. Let $i: K \hookrightarrow Y_{\infty}$, $1_K \mapsto \treeA$.
Observe that any tree $t$ can be witten as $t:=t_l \vee t_r= (t_l \succ \treeA) \nwarrow t_r$. Define recursively $\phi: \ Y_{\infty} \rightarrow LD $ by $\phi(\treeA):=x$ and for any tree $t\in Y_p$, $p>0$, $\phi(t):= (\phi(t_l) \succ x) \nwarrow \phi(t_r).$ Then, $\phi$ is the unique L-dipterous morphism extending $f$ and verifying $\phi \circ i=f$. Hence $Y_{\infty}$ is the free L-dipterous algebra over $K$.
\eproof

\noindent
As the L-dipterous operad is regular, we get the following corollary.
\begin{coro}
Let $V$ be a $K$-vector space. The $K$-vector space,
$$ \bigoplus_{n>0} \ KY_n \otimes V^{\otimes n},$$
equipped with the two following binary operations:
$$ (t \otimes \omega)\nwarrow (t' \otimes \omega'):= (t\nwarrow t') \otimes \omega\omega',$$
$$ (t \otimes \omega)\succ (t' \otimes \omega'):= (t\succ t') \otimes \omega\omega',$$
is the free L-dipterous algebra over $V$, the isomorphism $L-Dipt(V) \longrightarrow \bigoplus_{n>0} \ KY_n \otimes V^{\otimes n}$, being given by $v \mapsto \treeA \otimes v$.
\end{coro}

\Rk
We let the reader to find the free right L-dipterous algebra over a $K$-vector space $V$ by using the operation $t \nearrow s$ which consists in gluing the root of the tree $t$ on the most left leaf of $s$.

\Rk
There is an interesting link with duplicial-algebras. Such stuctures appear in the previous works of A. Brouder and A. Frabetti \cite{BF} and J.-L. Loday and M. Ronco \cite{LRbruhat}, see also \cite{Coa, Lerden} and \cite{Ltrip} for their extensions to the triplicial-algebras used in the good triple $(As, Trip, L)$ endowed with nonunital infinitesimal relations.
The free duplicial-algebras \cite{GB} is also constructed over rooted planar binary trees with the help of the operations $\nwarrow$ and $\nearrow$. A richer structure is then compatible with the underlying duplicial structure of planar rooted binary trees.
\subsection{From dendriform-Nijenhuis bialgebras}
Dendriform-Nijenhuis bialgebras have been introduced in \cite{NijLer} to produce  Baxter-Rota like operators, the so-called $TD$-operators, commuting with Nijenhuis operators.
A dendriform-Nijenhuis bialgebra is a triple $(A, \ \mu, \ \Delta )$ where $(A, \ \mu)$ is an associative algebra and $(A,\ \Delta )$ is a coassociative coalgebra such that,
$$ \Delta(ab) := \Delta(a)b + a\Delta(b) - \mu(\Delta(a)) \otimes b, \ \forall a,b \in A.$$
Set $x \bowtie y := \mu(\Delta(x)) y$ and $x \prec_{A} y := x \mu(\Delta(y))$ defined for all $x,y \in A$. Then the $K$-vector space $(A, \ \bowtie, \ \prec_{A} )$ is a right $L$-dipterous algebra.
\subsection{From preLie-algebras and quantum field theory}
In \cite{OudGuin}, Oudom and Guin construct explicitely over $Com^c(L)$ where $L$ is a PreLie-algebra two operations, $\circ$ and $*$, the last one being associative.
It is easy to observe that $(Com^c(L),\circ, *)$ is a right dipterous algebra, better the usual coproduct $\Delta$ of $Com^c$ is a morphism for these two operations.
This motivates the introduction of the triple of operads $(Com, Dipt, Prim_{Com} \ Dipt)$. We recall $(Com^c(L),\circ, *,\Delta)$ is dual to the Connes-Kreimer Hopf algebra introduced in quantum field theory.

\subsection{From rooted trees}
Only in this section, we consider labeled rooted trees (planarity is dropped, the degree of the nodes of trees can be equal to one and nodes are labeled). Rooted trees are registered under the name $A000169$ of the Encyclopedy of Integer Sequences and the number of labeled rooted trees with $n$ nodes is $n^{(n-1)}$.
For instance, $\bullet_v$, $\vert^w _v$ and so on.
The analogue of the grafting operation for planar rooted trees are here played by the so-called $NAP$-operation. $NAP$-algebras have been introduced via their coversions in \cite{Coa} under the name L-cocommutative coalgebra and independtly in \cite{Liv} where we borrow the terminology.
These are $K$-vector spaces equipped with a binary operation $\triangleleft$ verifying the following constraint:
$$ (x \triangleleft y ) \triangleleft z=(x \triangleleft z ) \triangleleft y. $$
On the $K$-vector space of rooted trees graduated by the number of nodes, define $t \triangleleft s$ to be the tree where the root of $s$ has been linked to the root of $t$. For instance $\bullet_v \triangleleft \bullet_w= \vert^w _v$. Then, any rooted tree $t$ whose root is labeled by $v$ can be decomposed in a unique way as:
$$ t:= (\ldots(v \triangleleft t_1)\triangleleft t_2) \triangleleft \ldots )\triangleleft t_n,$$
where the rooted trees $t_i$ are of smaller degrees. The free $NAP$-algebra over $V$ is then isomorphic to the $K$-vector space of rooted trees whose nodes are labeled by a basis of $V$ equipped with the operation $\triangleleft$ and
the canonical embedding,
$$i: V \hookrightarrow NAP(V), \ \ v \mapsto \bullet_v.$$
Recall also that a permutative algebra is a $K$-vector space \cite{Chp1} equipped with an associative product $\star$ verifying $x \star y \star z=x \star z \star y$. Consider the $K$-vector space $Perm(NAP(V))=NAP(V)\otimes (K\oplus Com(NAP(V)))$ and define the two following binary operations, $\star$ and $\prec$, for any trees $t, t_1, \ldots, t_n,r_1, \ldots, r_k,s, s_1 \ldots, s_p$ by:
$$ (t \otimes t_1 \ldots t_n) \star (s \otimes s_1 \ldots s_p):= t \otimes t_1 \ldots t_n s s_1 \ldots s_p, $$
%If $t:=(\ldots(\bullet_v \triangleleft t_1)\triangleleft t_2) \ldots )\triangleleft t_n)$, then,
$$ (t \otimes r_1 \ldots r_k) \prec (s \otimes s_1 \ldots s_p):=  (\ldots (t\triangleleft r_1) \triangleleft \ldots )\triangleleft r_k )\triangleleft s)\triangleleft \ldots ) \triangleleft s_p \otimes 1_K. $$
$$ (t \otimes 1_K) \star (s \otimes 1_K) =  t \otimes  s, $$
$$
(t \otimes 1_K) \star (s \otimes s_1 \ldots s_n )= t  \otimes ss_1 \ldots s_n,
$$
$$
(t \otimes t_1 \ldots t_n) \star (s \otimes 1_K) = t  \otimes t_1 \ldots t_n s,
$$
$$ (t \otimes 1_K) \prec (s \otimes s_1 \ldots s_p) = (\ldots(t\triangleleft  s) \triangleleft  s_1)\triangleleft \ldots ) \triangleleft s_p \otimes 1_K. $$
$$ (t \otimes t_1 \ldots t_n) \prec (s \otimes 1_K) = (\ldots(t \triangleleft  t_1)\triangleleft \ldots ) \triangleleft t_n ) \triangleleft s\otimes 1_K. $$
Because trees are non planar, the operation $\star$ is permutative, the operation $\prec$ is $NAP$ and the right dipterous axioms hold.

\section{The free dipterous algebra}
Let $V$ be a $K$-vector space.
The free dipterous algebra over $V$ was first announced in \cite{LodRon1}. We give here a proof as well as an explicit construction of the extension of linear maps into dipterous algebra morphisms.
Consider the $K$-vector space $\bar {T}(V)
:= \bigoplus_{n>0}\ V^{\otimes n}$.
Denote by $T_n$ the set of rooted planar trees with $n$ leaves. The cardinality of $T_n$ are registered under the name \textit{A001003 little Schroeder numbers} of the Online Encyclopedy of Integer Sequences. For $n=1,2,3$, we get:
$$ T_1=\{\ \vert \ \}, \ T_2=\{\ \treeA \ \}, \ T_3=\{\ \treeAB, \treeBA, \ \treeM \ \}.$$
Define grafting operations
by:
$$ T_{n_1} \times \ldots \times T_{n_p} \rightarrow T_{n_1 + \ldots + n_p}, \ \ (t_1, \ldots, t_p)\mapsto t_1\vee \ldots \vee t_p,$$
where the tree $t_1\vee \ldots \vee t_p$ is the tree whose roots of the $t_i$ have been glued together and a new root has been added. Observe that any rooted planar tree $t$ can be decomposed in a unique way via the grafting operation as
$t_{1} \vee \ldots \vee t_{p}$.
Set $T_{\infty}:= \bigoplus_{n>0} \ KT_n$.
Define over $\bar {T}(T_{\infty})$, the following binary operations, first on trees, then by bilinearity:
$$ (t_1 \otimes \ldots \otimes t_{n})\star (s_1 \otimes \ldots \otimes s_{p})=t_1 \otimes \ldots \otimes t_{n} \otimes s_1 \otimes \ldots \otimes  s_{p},$$
$$ (1) \ \ (t_1 \otimes \ldots \otimes t_{n})\succ s:= t_1 \vee  \ldots \vee t_{n} \vee  s_{1} \vee \ldots \vee s_{p}, \ \textrm{if,} \ s:=s_{1} \vee \ldots \vee s_{p},$$
$$(2)  \ (t_1 \otimes \ldots \otimes t_{n})\succ (s_1 \otimes \ldots \otimes s_{p})=t_1 \vee  \ldots \vee t_{n}\vee (s_1 \vee  \ldots \vee s_{p}), \ \textrm{if} \ p>1, $$
for all planar rooted trees $t_1, \ldots, t_n, s, s_1, \ldots, s_p$.
\begin{theo}
The free dipterous algebra over $K$ is the $K$-vector space $\bar{T}(T_{\infty})$ equipped with the usual concatenation operation $\star$ and $\succ$.
\end{theo}
\Proof
The concatenation product $\star$ is associative. To ease notation set $t_1 \otimes \ldots \otimes t_n:=t_1 \ldots t_n$, for any trees $t_1, \ldots, t_n$. Let $t_1, \ldots, t_n, s_1, \ldots, s_p, r_1, \ldots,  r_k$ be trees.\\
If $r:=r_1 \ldots  r_k$ with $k>1$, then:
\begin{eqnarray*}
(t_1 \ldots t_n \star s_1 \ldots s_p)\succ r &=& (t_1 \ldots t_n s_1 \ldots s_p)\succ r \\
&=& t_1 \vee  \ldots \vee t_{n}\vee s_1 \vee  \ldots \vee s_{p}\vee (r_1 \vee \ldots \vee r_k),
\end{eqnarray*}
and:
\begin{eqnarray*}
t_1 \ldots t_n \succ (s_1 \ldots s_p \succ r )&=& t_1 \ldots t_n \succ (s_1\vee \ldots\vee s_p \vee (r_1 \vee \ldots \vee r_k)) \ \textrm{(use (2)),}\\
&=& t_1 \vee  \ldots \vee t_{n}\vee s_1 \vee  \ldots \vee s_{p}\vee (r_1 \vee \ldots \vee r_k)\ \textrm{(use (1)).}
\end{eqnarray*}
If $r:=r_1 \vee \ldots \vee r_k$ with $k>0$, is a tree then:
$$t_1 \ldots t_n \succ (s_1 \ldots s_p \succ r )=t_1 \vee  \ldots \vee t_{n}\vee s_1 \vee  \ldots \vee s_{p}\vee r_1 \vee \ldots \vee r_k=(t_1 \ldots t_n \otimes s_1 \ldots s_p)\succ r.$$
Hence, $\bar {T}(T_{\infty})$ equipped with the concatenation operation $\star$ and $\succ$ is a dipterous algebra.\\ Let $t:=t_1 \vee \ldots \vee t_k$, $k>0$, be a tree with $t_k:= t_{k,1} \vee \ldots \vee t_{k,p}$. Observe that:
$$(3) \ \ t= t_1 \ldots t_{k-1} \succ (t_{k,1}  \ldots  t_{k,p}).$$
Let $(D, \star_D, \succ_D)$ be a dipterous algebra, $d \in D$ and $f:K \rightarrow D$, $1_K \mapsto d$.
Let $i: K \hookrightarrow \bar {T}(T_{\infty})$, $1_K \mapsto \vert$. Define recursively the map $\phi: \bar {T}(T_{\infty}) \rightarrow D$ by:
$$\phi(\vert)=d,$$
$$ \phi(t_1 \ldots t_n)=\phi(t_1)\star_D \ldots \star_D  \phi(t_n),$$
$$ (4) \ \ \phi(t)=\phi(t_1 \ldots t_{k-1} \succ (t_{k,1}  \ldots  t_{k,p}))= \phi(t_1 \ldots t_{k-1}) \succ_D \phi(t_{k,1}  \ldots  t_{k,p}).$$
Then, $\phi$ is a morphism of dipterous algebras. Indeed, this is by definition a morphism of associative algebras. Moreover, if $s_1 \ldots s_p$, $p>1$, is a forest, then:
\begin{eqnarray*}
\phi(t_1 \ldots t_n \succ s_1 \ldots s_p)&=&\phi(t_1 \vee \ldots \vee t_n \vee (s_1 \vee\ldots\vee  s_p)) \ \textrm{(use (2)),}\\
&=&\phi(t_1 \ldots t_n) \succ_D \phi(s_1 \ldots s_p) \ \textrm{(use (3-4)).}\\
\end{eqnarray*}
If $s$ is a tree, set $s:= s_1 \vee \ldots \vee s_p$, with $s_p:=s_{p,1} \vee \ldots \vee s_{p,k}$. Observe that
$\phi(\vert \succ \vert)=\phi(\vert) \succ_D \phi(\vert)$. We now proceed by induction.
\begin{eqnarray*}
\phi(t_1 \ldots t_n \succ s)&=&\phi(t_1\vee \ldots \vee t_n \vee s_1 \vee\ldots\vee  s_p)\ \textrm{(use (1)),}\\
&=&\phi(t_1 \ldots t_n s_1 \ldots s_{p-2} s_{p-1}\succ (s_{p,1}  \ldots  s_{p,k}))\ \textrm{(use (3)),}\\
&=&\phi(t_1 \ldots t_n s_1 \ldots s_{p-2} s_{p-1}) \succ_D \phi(s_{p,1}  \ldots  s_{p,k})\ \textrm{(use (4)),}\\
&=&(\phi(t_1 \ldots t_n s_1 \ldots s_{p-2})\star_D  \phi(s_{p-1})) \succ_D \phi(s_{p,1}  \ldots  s_{p,k}),\\
&=&\phi(t_1 \ldots t_n s_1 \ldots s_{p-2})\succ_D ( \phi(s_{p-1}) \succ_D \phi(s_{p,1}  \ldots  s_{p,k}))\ \textrm{(dipterous axioms),}\\
&=&\phi(t_1 \ldots t_n s_1 \ldots s_{p-2})\succ_D \phi(s_{p-1} \succ s_{p,1}  \ldots  s_{p,k}), \ \textrm{(by induction)} \ \\
&=&\phi(t_1 \ldots t_n s_1 \ldots s_{p-2})\succ_D \phi(s_{p-1} \succ s_{p,1}  \ldots  s_{p,k}),  \\
&=&\phi(t_1 \ldots t_n) \succ_D \phi(s_1\succ (s_2 \succ ( \ldots \succ( s_{p-2} \succ (s_{p-1} \succ (s_{p,1}  \ldots  s_{p,k}))\ldots),\\
& & \textrm{(by reapplying the process),}\\
&=&\phi(t_1 \ldots t_n) \succ_D \phi(s_1\vee s_2  \vee \ldots \vee  s_{p-1} \vee (s_{p,1}  \ldots  s_{p,k})),\\
& & \textrm{(by reapplying (1)),}\\
&=&\phi(t_1 \ldots t_n) \succ_D \phi(s_1\vee s_2  \vee \ldots \vee  s_{p-1} \vee s_{p} ) ,\\
&=&\phi(t_1 \ldots t_n) \succ_D \phi(s),
\end{eqnarray*}
proving by induction that $\phi$ is a morphism for the operation type $\succ$ and the only one such that $\phi\circ i=f$.
\eproof

\noindent
As the operad $Dipt$ is regular, the following holds.
Let $V$ be a $K$-vector space. The free dipterous algebra over $V$ is the $K$-vector space:
$$ Dipt(V):= \bigoplus_{n>0} Dipt_n \otimes V^{\otimes n}, $$
with  $Dipt(K):= \bigoplus_{n>0} Dipt_n\simeq \bar {T}(T_{\infty})$. Hence, $Dipt_n$ is explicitely described in terms of forests of rooted planar trees. By abuse of notation, we will mix $Dipt_n$ and its image under this isomorphism. Extend the operation $\star$ and $\succ$ as follows:
$$ ((t_1\ldots t_n) \otimes \omega) \star ((s_1\ldots s_p) \otimes \omega')= (t_1\ldots t_ns_1\ldots s_p) \otimes \omega \omega', $$
$$ ((t_1\ldots t_n) \otimes \omega) \succ ((s_1\ldots s_p) \otimes \omega')= (t_1\ldots t_n\succ s_1\ldots s_p) \otimes \omega \omega'.$$
The embedding map $i: V \hookrightarrow Dipt(V)$ is defined by: $v \mapsto \vert \otimes v$.

\noindent
Since the generating function associated with the Schur functor $\bar {T}$ is $f_{\bar {T}}(x):=\frac{x}{1-x}$ and with the Schur functor $T_{\infty}$ is
$f_{T_{\infty}}(x):=\frac{1 + x - \sqrt(1 - 6x + x^2)}{4}=x+x^2+3x^3+11x^4+45x^5 + \ldots$,
the generating function of the dipterous operad is $f_{\bar {T}}\circ f_{T_{\infty}}$, that is:
$$ f_{Dipt}(x):=\frac{1+x-\sqrt{1-6x+x^2}}{3-x+\sqrt{1-6x+x^2}}=x+2x^2+6x^3+22x^4+\ldots.$$
The sequence $(1,2,6,22,90,\ldots)$ is registered as $A006318 $ under the name \textit{Large Schroeder numbers}.

\NB
The free right dipterous algebra over a $K$-vector space $V$ is easily obtained from this construction (reverse the order in the definition of $\succ$).

\subsection{Rooted planar trees versus $m$-ary trees}
Fix $m>1$.
The aim of this subsection is to code rooted planar trees in terms of planar rooted $m$-ary trees in a injective way.
In \cite{Ler-polyg}, the notion of $m$-dendriform algebras, $m>1$ is introduced, the case $m=2$ being the dendriform algebras introduced by J.-L. Loday.
A $K$-vector space $T$ is a $m$-dendriform algebra if it is equipped with $m$ binary operations
$\prec,\succ,\bullet_2, \ldots, \bullet_{m-1}: T^{\otimes 2} \longrightarrow T$ verifying for all $x,y,z \in T,$ and for all $2 \leq i \leq m-1$,
the $\frac{m(m+1)}{2}$ axioms:
$$ (x\prec y)\prec z = x \prec (y \star z), \ \  \  \  \ (x \prec y)\bullet_i z =x\bullet_i (y \succ z) $$
$$ (x \succ y)\prec z = x \succ (y \prec z), \ \  \  \  \ (x \succ y)\bullet_i z =x\succ (y \bullet_i z) $$
$$(x \star y)\succ z =x\succ (y \succ z) , \ \  \  \  \ (x\bullet_i  y)\prec z = x \bullet_i (y \prec z), $$
where $x \star y := x\prec y + x\succ y$ and,
$$(x\bullet_i  y)\bullet_j z = x \bullet_i (y \bullet_j z), $$
for all $2 \leq i<j \leq m-1$.
The free $m$-dendriform algebra $m-Dend(V)$ over a $K$-vector space $V$ is related to planar rooted $m$-ary trees. As $m$-dendriform algebras are dipterous algebras we get:
$$\xymatrix{
V \ar[r]^-{i} \ar@{->}[rd]_{j} & Dipt(V)  \ar@{->}[d]^{\psi} \\
 & m-Dend(V)}$$
where $\psi$ is the unique morphism dipterous algebras  verifying $\psi\circ i=j$.
As a dipterous algebra is also a $m$-dendriform algebra by choosing all the other operations to be trivial, we get:
$$\xymatrix{
V \ar[r]^-{j} \ar@{->}[rd]_{i} & m-Dend(V) \ar@{->}[d]^{\Psi} \\
 & Dipt(V)}$$
with $\Psi$ the unique morphism of $m$-dendriform algebras such that $\Psi \circ j= i$. Therefore $(\Psi \circ \psi) \circ i=i$. As a morphism of $m$-dendriform algebras is also a morphism of dipterous algebras, we get by unicity: $\Psi \circ \psi=id_{_{Dipt(V)}}$.

\section{On the good triple of operads $(As, Dipt, Mag^{\infty} )$}
The aim of this section is to provide a Poincar\'e-Birkhoff-Witt like theorem and a Cartier-Milnor-Moore like theorem for infinitesimal dipterous bialgebras.
By definition, an infinitesimal dipterous bialgebra $(D, \succ,\star, \Delta)$ is a dipterous algebra equipped with a coassociative coproduct $\Delta$ verifying the following so-called nonunital semi-infinitesimal relations:
$$ \Delta(x \succ y):= x_{(1)}\otimes (x_{(2)} \succ y) + (x \star y_{(1)}) \otimes y_{(2)} + x\otimes y.$$
$$ \Delta(x \star y):= x_{(1)}\otimes (x_{(2)} \star y) + (x \star y_{(1)}) \otimes y_{(2)} + x\otimes y.$$
It is said to be connected when $D=\bigcup_{r>0} F_rD$ with
the filtration $(F_rD)_{r>0}$ defined as follows:
$$ \textrm{(The primitive elements)} \ \ F_1D:=Prim \ D=\{x \in D, \ \Delta(x)=0 \},$$
Set $\Delta^{(1)}:= \Delta$ and
$\Delta^{(n)}:= (\Delta \otimes id_{n-1})\Delta^{(n-1)}$ with
$id_{n-1}= \underbrace{id \otimes \ldots \otimes id}_{times \ n-1}$. Then,
$$ F_rD:=\ker \ \Delta^{(r)}.$$

\begin{theo}
\label{h1}
Let $V$ be a $K$-vector space. Define on $ Dipt(V)$,
the free dipterous algebra over $V$,
the cooperation $\Delta: Dipt(V) \rightarrow Dipt(V) \otimes Dipt(V)$ recursively as follows:
$$ \Delta(i(v)):=0, \ \textrm{for all} \ v \in V,$$
$$ \Delta(x \succ y):= x_{(1)}\otimes (x_{(2)} \succ y) + (x \star y_{(1)}) \otimes y_{(2)} + x\otimes y.$$
$$ \Delta(x \star y):= x_{(1)}\otimes (x_{(2)} \star y) + (x \star y_{(1)}) \otimes y_{(2)} + x\otimes y,$$
for all $x,y \in Dipt(V)$.
Then $(Dipt(V), \Delta)$
is a connected infinitesimal  dipterous bialgebra.
\end{theo}
\Proof
We show by induction that $\Delta$ is coassociative. Coassociativity holds over $K\vert \otimes V$. As $Dipt$ is a binary operad, if $z \in Dipt(V)$, then there exists $x,y \in Dipt(V)$ with smaller degrees and such that $z=x \diamond y$, where $\diamond= \succ; \ \star$ and:
\begin{eqnarray*}
\Delta(x \diamond y)= x_{(1)}\otimes (x_{(2)} \diamond y)+ (x \star y_{(1)}) \otimes y_{(2)} + x \otimes y.
\end{eqnarray*}
Therefore,
\begin{eqnarray*}
(\Delta \otimes id)\Delta(x \diamond y)&= & x_{(11)}\otimes x_{(12)} \otimes(x_{(2)} \diamond y)+
x_{(1)}\otimes (x_{(2)} \star y_{(1)}) \otimes y_{(2)}
+ (x \star y_{(11)}) \otimes y_{(12)} \otimes y_{(2)}\\
& & +
x \otimes y_{(1)} \otimes y_{(2)}
+ x_{(1)}\otimes x_{(2)}  \otimes y,
\end{eqnarray*}
and,
\begin{eqnarray*}
(id \otimes \Delta)\Delta(x \diamond y)&=&
x_{(1)}\otimes x_{(21)} \otimes (x_{(22)}\diamond y)+
x_{(1)}\otimes (x_{(2)} \star y_{(1)}) \otimes y_{(2)}+
x_{(1)}\otimes x_{(2)} \otimes y \\
& & +
(x \star y_{(1)}) \otimes y_{(21)}\otimes y_{(22)}
+x \otimes y_{(1)} \otimes y_{(2)}.
\end{eqnarray*}
Hence the coassociativity of $\Delta$ by induction.
We have to show that the coproduct respects the dipterous axioms:
\begin{eqnarray*}
\Delta(x \diamond (y\diamond z))&=&
x_{(1)} \otimes (x_{(2)} \diamond (y\diamond z)) +
(x \star y_{(1)}) \otimes (y_{(2)}\diamond z)\\
& & +
(x \star y \star z_{(1)}) \otimes  z_{(2)} +
(x \star y) \otimes z + x \otimes (y \diamond z),
\end{eqnarray*}
and,
\begin{eqnarray*}
\Delta((x \star y)\diamond z)&=&
x_{(1)} \otimes ((x_{(2)}\star y)\diamond z)+
(x \star y_{(1)}) \otimes (y_{(2)}\diamond z)+
x\otimes (y\diamond z)\\
& & +
(x \star y \star z_{(1)}) \otimes z_{(2)}
+
(x \star y)\otimes z.
\end{eqnarray*}
Hence, the coproduct is compatible with the dipterous axioms. Connectedness follows by construction, hence the result.
\eproof

\NB
For the right dipterous case $(AD, \star, \prec)$, we have to consider the following compatibility relation:
$$ \Delta(x \prec y):= x_{(1)}\otimes (x_{(2)} \star y) + (x \prec y_{(1)}) \otimes y_{(2)} + x\otimes y.$$
$$ \Delta(x \star y):= x_{(1)}\otimes (x_{(2)} \star y) + (x \star y_{(1)}) \otimes y_{(2)} + x\otimes y,$$
and all our results will still hold.

\Rk
Following the works \cite{Rota, AguiarLoday, Lertribax}, one can enlarge the definition of the nonunital semi-infinitesimal relation to include a parameter $t \in K$ as follows for instance for the right dipterous case,
$ \Delta(x \prec y):= x_{(1)}\otimes (x_{(2)} \star y) + (x \prec y_{(1)}) \otimes y_{(2)} + tx\otimes y$ and
$\Delta(x \star y):= x_{(1)}\otimes (x_{(2)} \star y) + (x \star y_{(1)}) \otimes y_{(2)} + tx\otimes y.$

\noindent
We now apply Theorem 2.5.1 of Loday \cite{GB} and his notation to recover a P.B.W. and C.M.M. like theorems for infinitesimal dipterous bialgebras.
The hypothesis $(H0)$ is verified since the semi-nonunital infinitesimal relations are distributive. The hypothesis $(H1)$ also holds because of Theorem \ref{h1}. Then, Theorem 2.2.2 \cite{GB} claims that the Schur functor $\mathcal{P}$ given by $\mathcal{P}(V) :=(Prim_{As} \ Dipt)(V)$ is a suboperad of the operad $Dipt$. Hence, the forgetful functor $F: Dipt-alg \rightarrow \mathcal{P}-alg$ defined in Section 2.4.3 \cite{GB}
has a left adjoint $U: \mathcal{P}-alg \rightarrow  Dipt-alg$, called the universal enveloping algebra functor.
\begin{theo}
\label{isom}
For any infinitesimal dipterous bialgebras $D$, the following are equivalent:
\begin{enumerate}
\item{$D$ is connected;}
\item{$D$ is isomorphic to $U(Prim \ D)$ as an infinitesimal dipterous bialgebra;}
\item{$D$ is cofree among connected coassociative coalgebras: $D \simeq As^c(Prim \ D)$. }
\end{enumerate}
\end{theo}
\Proof
We have to check hypothesis $(H2epi)$ of Theorem 2.5.1 \cite{GB}. Let $V$ be a $K$-vector space. The cofree coassociative coalgebra over $V$ is $As^c(V)= \bigoplus_{n>0} \ V^{\otimes n}$ as a $K$-vector space equipped with the deconcatenation coproduct
$\delta$ given by:
$$ \delta(v_1 \otimes \ldots \otimes  v_n):= \sum_{k=1}^{n-1} \ (v_1 \otimes \ldots \otimes  v_k) \otimes (v_{k+1} \otimes \ldots \otimes v_n),$$
for all $v_1 \ldots v_n \in V$.
Therefore, there exits a unique coalgebra morphism $\phi(V): Dipt(V) \rightarrow As^c(V)$ such that $\pi \circ\phi(V)=\pi',$ where $\pi: As^c(V) \twoheadrightarrow V$ and $\pi':  Dipt(V) \twoheadrightarrow V$ are the canonical projections. This morphism is surjective since it maps any trees of $\bar {T}(T_{\infty})$ into $1_K$. The map: $$s(V):As^c(V) \rightarrow Dipt(V) $$
$$(v_1 \otimes \ldots \otimes  v_n) \mapsto  \underbrace{(\vert \vert \ldots \vert)}_{times \ n}\otimes (v_1 \otimes \ldots \otimes  v_n), $$
is obviously a coalgebra morphism and verifies $\phi(V) \circ s(V)=id$. Hence the three hypotheses of Theorem 2.5.1 \cite{GB} are checked, hence the results.
\eproof
\NB
Similarly, one can obtain a unital version of this result, see Section 7.
\section{The operad $Mag^{\infty}$}
\noindent
We now explicit the primitive operad $\mathcal{P}$.
Recall a $Mag^{\infty}$-algebra $G$ is a $K$-vector space equipped with one $n$-ary operation for each $n>1$, denoted by:
$$[ \ , \ldots, \ ]_n: G^{\otimes n} \rightarrow G.$$

\begin{prop}
The $K$-vector space $T_{\infty}$
equipped with the grafting operations is the free $Mag^{\infty}$-algebra over $K$.
The generating function of the operad $Mag^{\infty}$ is $f_{Mag^{\infty}}(x):=f_{T_\infty}(x)$.
\end{prop}
\Proof
The operad $Mag^{\infty}$ is regular.
Consider now the map $Mag^{\infty}(K) \rightarrow T_{\infty}$, $[ \ , \ldots, \ ]_n \mapsto cor_n,$
where $\treeM_n$ denotes the $n^{th}$ corolla. For instance, $cor_2=\treeA$, $cor_3=\treeM$, $cor_4=\treeCor$ and so forth. This map is an isomorphism of $Mag^{\infty}$-algebras, hence the claim.
\eproof

\begin{prop}
\label{grovve}
The primitive part of a connected infinitesimal dipterous bialgebra is a $Mag^{\infty}$-algebra.
\end{prop}
\Proof
Set $\sqsubset x,y \sqsupset_2= x \succ y- x \star y:= x \triangleleft y$, and
$$\sqsubset x_1, \ldots, x_n\sqsupset_{n}= x_1 \triangleleft  (x_2 \succ (x_3 \succ  \ldots \succ(x_{(n-2)}\succ (x_{(n-1)}\succ x_n)\ldots),$$
for all $n>0$ and $x, x_i \in D$. If $x_1, \ldots, x_n  \in Prim D$ so is $\sqsubset x_1, \ldots, x_n \sqsupset_{n}$ because,
\begin{eqnarray*}
\Delta(x_1 \succ  (x_2 \succ  \ldots \succ (x_{(n-1)}\succ x_n)\ldots)&=& x_1\star \Delta (x_2 \succ (x_3 \succ  \ldots \succ (x_{(n-1)}\succ x_n)\ldots)\\
& & + x_1\otimes (x_2 \succ (x_3 \succ  \ldots \succ (x_{(n-1)}\succ x_n)\ldots)\\
&=& \Delta(x_1 \star  (x_2 \succ  \ldots \succ (x_{(n-1)}\succ x_n)\ldots).
\end{eqnarray*}
Therefore, $Prim_{As} \ Dipt$-algebras are $Mag^{\infty}$-algebras. Because of Theorem \ref{isom}, there is an isomorphism of Schur functors,
$$ Dipt = As^c \circ Prim_{As} \ Dipt, $$
Hence, we have $\dim Prim_{As} \ Dipt(n)=\dim Mag^{\infty}(n)$, for all $n>0$. However $Dipt=As \circ Mag^\infty$ by construction. Since $(As,As, Vect)$ endowed with the infinitesimal relation is good, $Prim \ As(Mag^\infty(V))=Mag^\infty(V)$. Hence,
the operad of primitive elements $Prim_{As} \ Dipt$ is the operad $Mag^{\infty}$.
\eproof

\NB
We can now explicit the universal enveloping algebra functor $U: Mag^{\infty}-alg. \rightarrow Dipt-alg.$. Let $G$ be a $Mag^{\infty}$-algebra whose $n$-ary operations for all $n>1$ are denoted by
$[ \ , \ldots, \ ]_n: G^{\otimes n} \rightarrow G.$ Then, $U(G)$ is the quotient of $Dipt(G)$ by the relations which consist in identifying the operations $[ \ , \ldots, \ ]_n,$ $n>1,$ to the operations $\sqsubset \ , \ldots, \  \sqsupset_n$ of $Dipt(G)$ made on $\star$ and $\succ$.

\NB Let $D$ be a connected infinitesimal dipterous bialgebra. Then, $e: D \rightarrow Prim \ D$ defined recursively by $x \mapsto e(x):= x-x_{(1)} \star e(x_{(2)})$ is an idempotent. (Indeed, the proof of Proposition 2.5 \cite{LodRon} holds in our case since $\star$ is associative.)

\noindent
According to the terminology developed in \cite{GB}, we summarize our results in the following theorem.
\begin{theo}
The triple of operads $(As, Dipt,Mag^{\infty} )$ endowed with the nonunital semi-infinitesimal relations is good.
\end{theo}

\section{More on dipterous algebras}
The aim of this section is to give an homology theory for dipterous algebras and to prove that the operad $Dipt$ is Koszul.
\subsection{The quasi-nilpotent dipterous operad}
The dual in the sense of Ginzburg and Kapranov \cite{GK} of the operad $Dipt$ is $Dipt^!:=QNDipt$, the so-called quasi-nilpotent dipterous operad.
A $QNDipt$-algebra $Q$ is a $K$-vector space equipped with two binary operations $\star$ and $\succ$ verifying dipterous axioms plus:
\begin{eqnarray}
(x\succ y)\star z &=& 0,\\
x\succ (y\star z) &=& 0,\\
x \star(y\succ z) &=& 0,\\
(x\succ y)\succ z &=& 0,
\end{eqnarray}
for all $x,y,z \in Q$.
\begin{theo}
Let $V$ be a $K$-vector space. Then,
the $K$-vector space,
$$ QNDipt(V):=T(V)\otimes (K \oplus V),$$
equipped with the following two operations, defined by:
$$ X\otimes 1_K \star Y \otimes 1_K := XY\otimes 1_K,$$
$$ X\otimes 1_K \succ Y \otimes v := XY\otimes v, \ \textrm{if} \ Y \notin V,$$
$$ X\otimes 1_K \succ v \otimes 1_K := X\otimes v,$$
for all $X,Y \in T(V)$ and $v \in V$, vanishing otherwise, and
then extended by bilinearity is the free $QNDipt$-algebra over $V$. Moreover, the generating function of the operad $QNDipt$ is:
$$ f_{QNDipt}(x):=x \ \frac{1+x}{1-x}=x + 2x^2 + 2x^3 +\ldots+2x^n +\ldots .$$
\end{theo}
\Proof
So defined, the operations $\star$ and $\succ$ obey $QNDipt$ axioms. For instance, we check Axiom $(3)$:
\begin{eqnarray*}
X\otimes x \succ (Y \otimes y \star Z \otimes z) &=& X\otimes x \succ (YZ \otimes 1_K) \ (Vanish  \ except \ if \ y=z=1_K)\\
&=& 0.
\end{eqnarray*}
Set $i: V \hookrightarrow QNDipt(V)$, $v \mapsto v\otimes 1_K$.
Let $Q$ be a $QNDipt$-algebra and $f:V \mapsto Q$ be a linear map. Let $v_1, \ldots, v_n,v \in V$. Then,
$$ v_1 \ldots v_n \otimes 1_K = i(v_1)\star \ldots \star i(v_n), $$
$$v_1 \ldots v_n \otimes v= (i(v_1)\star \ldots \star i(v_n))\succ i(v).$$
Define the map $\phi: QNDipt(V) \rightarrow Q$ by $\phi(i(v))=f(v)$ and by,
$$ \phi(v_1 \ldots v_n \otimes 1_K):= f(v_1) \star_Q \ldots \star_Q f(v_n),$$
$$\phi(v_1 \ldots v_n \otimes v):=
(f(v_1) \star_Q \ldots \star_Q f(v_n))\succ_Q f(v).$$
Then, $\phi$ is the isomorphism of $QNDipt$-algebras extending $f$.
Checking the morphism of associative algebras property is left to the reader. Let $X,Y \in T(V)$ and $x,y \in V$ or $K$.
set $Y=y_1 \ldots y_n$, $n>1$.
On the one hand,
$\phi(X\otimes x \succ Y \otimes y)=
0$ unless $x=1_K$ and $y\not= 1_K$. In this case, $\phi(X\otimes 1_K \succ Y \otimes y)= \phi(X Y \otimes y)=(\phi(X\otimes 1_K)\star_Q \phi(Y\otimes 1_K))\succ_Q f(y)$. On the other hand,
$\phi(X\otimes x) \succ_Q \phi(Y \otimes y)= (\phi(X\otimes 1_K) \succ_Q  f(x)) \succ_Q (\phi(Y\otimes 1_K) \succ_Q f(y)),$ if $x\not=1_K$ and $y\not=1_K$. But this vanishes because of Axiom (4). If $x=1_K$ and $y\not=1_K$, then
$\phi(X\otimes 1_K) \succ_Q \phi(Y \otimes y)= \phi(X\otimes 1_K) \succ_Q (\phi(Y\otimes 1_K) \succ_Q f(y))=(\phi(X\otimes 1_K) \star_Q \phi(X\otimes 1_K)) \succ_Q f(y).$
For the last case suppose
 $y=1_K$, then
$\phi(X\otimes 1_K) \succ_Q \phi(Y\otimes 1_K)= \phi(X\otimes 1_K) \succ_Q ((f(y_1)\star_D \ldots \star_D f(y_{(n-1)})\star_D f(y_n)) =0$ because of Axiom (2). Hence, the required equality if $n>1$.

For $Y\in V$, the same proof gives
that $\phi(X\otimes x) \succ_Q \phi(Y \otimes y)$ vanishes except if $x=y=1_K$. In this case,
$\phi(X\otimes 1_K) \succ_Q \phi(Y \otimes 1_K)=\phi(X\otimes 1_K) \succ_Q f(Y).$ Similarly,
$\phi(X\otimes x \succ_Q Y \otimes y)$ vanishes except if $x=y=1_K$ since $Y \in V$. In this case,
$\phi(X\otimes x \succ_Q Y \otimes y)=\phi(X\otimes Y)=\phi(X\otimes 1_K) \succ_Q f(Y)$, hence the linear map $\phi$ is a morphism of
$QNDipt$-algebras, the only one extending $f$.
The last claim is straightforward since the generating function associated to the Schur functor $\bar {T}$ is $f_{\bar {T}}(x):=\frac{x}{1-x}$.
\eproof

\subsection{Homology of dipterous algebras}
Because $QNDipt$ is regular we have $QNDipt(n):=QNDipt_n \otimes K\Sigma_n$, with $\dim QNDipt_n= 2$ for $n>1$ and $\dim QNDipt_1= 1$. Let $\mathcal{P}$ be an operad. According to \cite{GK}, a chain complex can be constructed for any $\mathcal{P}$-algebras with the help of the dual operad $\mathcal{P}^!$.
In our case, let $(D, \star, \succ)$ be a dipterous algebra. Set $S:=\{\divideontimes, \rhd \}$.
Define the module of $n$-chains, $n>1$,  by:
$$ C_n(D):= KS \otimes D^{\otimes n}$$
$$ C_1(D):=D.$$
Fix $n>1$. For any $1 \leq i\leq n-1$, define linear maps (face maps)  $d_i: C_n(D) \rightarrow C_{n-1}(D)$ by:
$$ d_i(\divideontimes \otimes x_1 \ldots x_n):= \divideontimes \otimes x_1 \ldots x_{(i-1)} (x_i \star x_{i+1})x_{(i+2)} \ldots x_n,$$
for $1 \leq i \leq n-1$; By:
$$ d_i(\rhd \otimes x_1 \ldots x_n):= \rhd  \otimes x_1 \ldots x_{(i-1)} (x_i \star x_{i+1})x_{(i+2)} \ldots x_n,$$
for $1 \leq i < n-1$
and:
$$ d_{n-1}(\rhd \otimes x_1 \ldots x_n):= \rhd \otimes x_1 \ldots x_{n-2} (x_{n-1} \succ x_n).$$
Set for all $n>1$:
$$d: = \sum_{i=1}^{n-1} \ (-1)^{i+1} d_i: C_n(D) \rightarrow C_{n-1}(D).$$
\begin{prop}
For all $n>1$,
the face maps $d_i: C_n(D) \rightarrow C_{n-1}(D)$ obey the simplicial relations:
$$ d_i d_j= d_{j-1} d_i, \ \ for \ any \ 1\leq i < j < n.$$
Moreover, $d^2=0$ and,
$$ (C_*(D),d):\ \  \cdots \xrightarrow{d} KS \otimes D^{\otimes n} \xrightarrow{d} KS \otimes D^{\otimes n-1} \xrightarrow{d} \cdots \xrightarrow{d} KS \otimes D^{\otimes 2} \xrightarrow{d} D,$$
is a chain-complex.
\end{prop}
\Proof
Fix $n>2$. Since $\star$ is associative, the face maps restricted to the symbol $\divideontimes$ obey the simplicial relations (Hochschild Homology). Let $1\leq i < j < n$. Restricted to the symbol $\rhd$, the simplicial relations hold if $j> i+1$ or if $1\leq i<j < n-1$. For $i=n-1$, we get:
$$ \rhd \otimes x_1 \ldots x_n \xrightarrow{d_n} \rhd \otimes x_1 \ldots x_{n-2}(x_{n-1} \succ x_n) \xrightarrow{d_{n-1}} \rhd \otimes x_1 \ldots x_{n-3} (x_{n-2}\succ(x_{n-1} \succ x_n)), $$
$$ \rhd \otimes x_1 \ldots x_n \xrightarrow{d_{n-1}} \rhd \otimes x_1 \ldots (x_{n-2}\star x_{n-1}) x_n \xrightarrow{d_{n-1}} \rhd \otimes x_1 \ldots x_{n-3} ((x_{n-2}\star x_{n-1})\succ x_n), $$
hence $d_{n-1}d_n=d_{n-1}d_{n-1}$ and the simplicial relations hold.
\eproof

\noindent
By definition, the homology of a dipterous algebra $D$ is the homology of the chain-complex $(C_*(D),d)$:
$$  HC_n(D):=H_n(C_*(D),d), \ \ n \geq 1.$$
By definition, the cohomology of a dipterous algebra $D$ is:
$$  HC^n(D):=H^n(Hom(C_*(D),K)), \ \ n \geq 1.$$
\begin{theo}
Let $V$ be a $K$-vector space. Then,
$H_1(Dipt(V)) \simeq V$ and $H_n(Dipt(V))=0$ for $n>1$, that is
the operad $Dipt$ is Koszul.
\end{theo}
\Proof
By definition $H_1(Dipt(V)) = Dipt(V)/\{ x \star y, \ x \succ y \ | \ x,y \in Dipt(V) \} \simeq V.$
To show that the complex above is acyclic, define the following map,
$$ h=h_n:KS \otimes Dipt(V)^{\otimes n}\rightarrow KS \otimes Dipt(V)^{\otimes (n+1)},$$
by $h_n(\rhd \otimes x_1 \ldots x_{n-1}u):= (-1)^{n+1} \ \rhd \otimes x_1 \ldots x_{n-1} A b$
if $u$ can be written in a unique way as,
$$ u:=a_p \succ (a_{p-1}\succ( \ldots \succ (a_2 \succ (a_1 \succ b))\ldots),$$
(always possible since $Dipt(V)$ is the free dipterous algebra over $V$) for some $a_1, \ldots, a_p,b \in Dipt(V)$ and where $A:=a_p \star  \ldots \star a_1$ and vanishes otherwise and
$h_n(\divideontimes \otimes x_1 \ldots x_{n-1}u):= (-1)^{n+1}  \ \divideontimes \otimes x_1 \ldots Ab,$
if $u$ can be written in a unique way as,
$$ u:= a_p \star \ldots \star a_1 \star b,$$
for some $a_1, \ldots, a_p,b \in Dipt(V)$,
and where $A:=a_p \star \ldots \star a_1$
and vanishes otherwise.
Then $h$ is an homotopy and the relation $dh_n + h_{n-1}d=id$ holds.
Hence $H_n(Dipt(V))=0$ for $n>1$. The fact the operad $Dipt$ is Koszul is then a direct consequence of \cite{GK}.
\eproof

\section{The good triple $(2As, Dipt, Vect)$}
The aim of this section is to prove that the triple of operads $(2As, Dipt, Vect)$ endowed both with the unital semi-Hopf and with the unital semi-infinitesimal compatibility relations explained below is good.
First we follow \cite{LodRon}.
A unital dipterous algebra, i.e., a dipterous algebra $(D, \star, \succ)$ equipped with a unit $1$ verifying
the Loday-Ronco axioms: $1\star x=x=x\star 1$ and $1\succ x=x$ and $x \succ 1=0$ for all $x \in D$, the symbol $1 \succ 1$ being not defined. For instance $K \oplus Dipt(V)$, with $V$ a $K$-vector space, is a unital dipterous algebra with unit $1_K$.
If $A$ and $B$ are two unital dipterous algebras then one can define a unital dipterous algebra structure over $A \otimes B$ as follows:
\begin{eqnarray*}
(a \otimes b) \star (a' \otimes b') &:=&(a \star a') \otimes (b \star b'),\\
(a \otimes b) \succ (a' \otimes b') &:=&(a \star a') \otimes (b \succ b'), \ \textrm{if} \ b \otimes b' \not= 1_A \otimes 1_B, \\
(a \otimes 1_A) \succ (a' \otimes 1_B) &:=& (a \succ a') \otimes 1_B, \
\end{eqnarray*}
for all $a,a' \in A$ and $b,b' \in B$.

\NB
Recall that the classical structure requires:
$$(a \otimes b) \diamond (a' \otimes b') :=(a \diamond a') \otimes (b \diamond b'),$$
for any generating operation $\diamond$.
Hence, the classical structure coincides with this structure only for associative operations. We use the word ``\textsf{semi}'' to refer to its unusual structure.

Recall the operad $2As$ has been introduced in \cite{LodRon}. A $K$-vector space equipped with two associative products is called a 2-associative algebra. It has been shown in \cite{LodRon} that the free 2-associative algebra over a $K$-vector space $V$ is related also to rooted planar trees. In fact, we have $\dim \ 2As_n=\dim \ Dipt_n$ for all $n>0$.

By mixing our results in Section 4 and results of \cite{LodRon}, we get another interesting notion of bialgebras
we name $2As^c-Dipt$-bialgebras. Such an object is:
\begin{enumerate}
\item{A unital dipterous algebra $(\mathcal{H}, \star, \succ)$.}
\item{Plus two coassociative coproducts $\vartriangle, \ \blacktriangle: \mathcal{H}\rightarrow \mathcal{H}^{\otimes 2}$ verifying the following compatibility relations for all $x,y \in \mathcal{H}$:
\begin{enumerate}
\item{The so-called unital semi-infinitesimal relations for $\vartriangle$:
$$ \vartriangle(1)=1 \otimes 1.$$
$$ \vartriangle(x \succ y):=  \vartriangle(x) \succ (1 \otimes y) + (x \otimes 1) \succ \vartriangle(y)  - x\otimes y.$$
$$ \vartriangle(x \star y):= \vartriangle(x) \star (1\otimes y) + (x \otimes 1)  \star \vartriangle(y)  - x\otimes y.$$
}
\item{The unital semi-Hopf relations for $\blacktriangle$ \cite{LodRon}:
$$ \blacktriangle(1)=1 \otimes 1,$$
$$\blacktriangle(x \star y)= \blacktriangle(x) \star \blacktriangle(y),$$
$$\blacktriangle(x \succ y)= \blacktriangle(x) \succ \blacktriangle(y).$$
}
\end{enumerate}

}
\end{enumerate}

We now follow \cite{GB}. Such a bialgebra $\mathcal{H}$ is said to be connected if $\mathcal{H}= \cup_{r \geq 0} \ \mathcal{F}_r$ where the filtration $(\mathcal{F}_r)_{r \geq 0}$ is defined as follows. First introduced the reduced coproducts:
$$ \Delta := \vartriangle- 1_K \otimes id - id \otimes 1_K,$$
$$ \bar{\blacktriangle} := \blacktriangle - 1_K \otimes id - id \otimes 1_K.$$
Then,
$$ \mathcal{F}_1:=Prim \ \mathcal{H}:=\{x \in \mathcal{H}; \ \Delta(x)=0=\bar{\blacktriangle}(x)\},$$
$$ \mathcal{F}_r:=\{x \in \mathcal{H}; \ \forall \ n\geq r, \ \Delta^{(n)}(x)=0=\bar{\blacktriangle}^{(n)}(x)\}.$$
We add $\mathcal{F}_0:=K.1_{\mathcal{H}}.$
Observe that such a connected $2As^c-Dipt$-bialgebra has a usual counit $\epsilon: \mathcal{H}\rightarrow K$ defined by $\epsilon(1_{\mathcal{H}})=1_K$ and $\epsilon(x)=0$ for all $x \in \mathcal{H}$ different from $1_{\mathcal{H}}$.
\begin{prop}
\label{2as-dipt}
Let $V$ be a $K$-vector space. Then the unital free dipterous algebra over $V$ is a connected $2As^c-Dipt$-bialgebra.
\end{prop}
\Proof
The $K$-vector space $Dipt(V)\oplus K$ is the unital free dipterous algebra over $V$ with unit $1_K$ and with embedding $i: V\hookrightarrow Dipt(V)$.

For all $v \in V$ set:
$$ \blacktriangle(i(v)):= 1_K \otimes i(v) + i(v) \otimes 1_K.$$
$Dipt(V)\oplus K$ is free, there exists a unique dipterous algebra extension $\blacktriangle:  Dipt(V) \rightarrow (Dipt(V)\oplus K)^{\otimes 2}$  of the map $ V \rightarrow (Dipt(V)\oplus K)^{\otimes 2}, \ v \mapsto 1_K \otimes i(v) + i(v) \otimes 1_K$
\cite{LodRon}. Add now $\blacktriangle(1_K):=1_K \otimes 1_K$ to get the first coassociative coproduct.
For the other one, set
for all $v \in V$:
$$ \vartriangle(1_K)=1_K \otimes 1_K,$$
$$ \vartriangle(i(v)):= 1_K \otimes i(v) + i(v) \otimes 1_K,$$
$$ \vartriangle(x \succ y):=  \vartriangle(x) \succ (1 \otimes y) + (x \otimes 1) \succ \vartriangle(y)  - x\otimes y,$$
$$ \vartriangle(x \star y):= \vartriangle(x) \star (1\otimes y) + (x \otimes 1)  \star \vartriangle(y)  - x\otimes y,$$
for all
$x,y \in Dipt(V)$.
For $x,y \in Dipt(V)$, using the action of the unit $1_K$ on $(K \oplus Dipt(V))^{\otimes 2}$ we get for instance:
\begin{eqnarray*}
\vartriangle(x \succ y) &=& (1_K \otimes x + x \otimes 1_K + x_{(1)} \otimes x_{(2)})\succ (1_K \otimes y) +  (x \otimes 1_K)\succ (1_K \otimes y + y \otimes 1_K \\
& & + y_{(1)} \otimes y_{(2)})- x\otimes y,\\
&=&  1_K \otimes (x \succ y) + (x \star 1_K) \otimes (1_K \succ y) + x_{(1)} \otimes (x_{(2)}\succ y) + (x \star 1_K) \otimes (1_K \succ y) \\ & &
+ (x \succ y) \otimes 1_K
 + (x \star y_{(1)}) \otimes y_{(2)}-x \otimes y,\\
&=& (x \succ y) \otimes 1_K + 1_K \otimes (x \succ y)+ [x_{(1)} \otimes (x_{(2)}\succ y)+ (x \star y_{(1)}) \otimes y_{(2)} + x \otimes y].
\end{eqnarray*}
The reduced coproducts:
$$ \Delta := \vartriangle- 1_K \otimes id - id \otimes 1_K,$$
$$ \bar{\blacktriangle} := \blacktriangle - 1_K \otimes id - id \otimes 1_K,$$
can be introduced. Observe then that the coproduct $\Delta$ is the one introduced in Section 4. The fact that
$\Delta(1_K \succ x)=\Delta(x)$ and $\Delta(x \succ 1_K)=0$
is straightforward and show that $\Delta$ well-behaves with the action of the unit $1_K$.
One then define as expected the filtration $(\mathcal{F}_r)_{r \geq 0}$ as follows:
$$ \mathcal{F}_0:=K.1_K,$$
$$ \mathcal{F}_1:=Prim \ \mathcal{H}:=\{x \in \mathcal{H}; \ \Delta(x)=0=\bar{\blacktriangle}(x)\},$$
$$ \mathcal{F}_r:=\{x \in \mathcal{H}; \ \forall \ n\geq r, \ \Delta^{(n)}(x)=0=\bar{\blacktriangle}^{(n)}(x)\}.$$
By construction, $(K \oplus Dipt(V))$ is a
connected $2As^c-Dipt$-bialgebra.
\eproof

\noindent
We now adapt a Quillen's result to our $2As^c-Dipt$-bialgebras.
\begin{lemm}
\label{Quillen}
Let $\mathcal{H}, \mathcal{H}'$ be  connected $2As^c-Dipt$-bialgebras and $\theta: \mathcal{H} \rightarrow \mathcal{H}'$ be a morphism of graduated $2As^c$-coalgebras. Then $\theta$ is injective on $\mathcal{H}$ if and only if
it is injective on $Prim \ \mathcal{H}$.
\end{lemm}
\Proof
Observe that $\theta(1_\mathcal{H})=1_\mathcal{H'}$.
We suppose $\theta$ to be injective on $Prim \ \mathcal{H}:=\mathcal{F}_1\mathcal{H}$ and not injective on $\mathcal{H}$. Let $x \in  \ker \theta \cap \mathcal{F}_{r}\mathcal{H} $ different from zero plus  with $r$ minimal. Thus $r>1$. Using one reduced coproduct we can write for instance $(\theta\otimes \theta)\Delta(x)=\theta(x_{(1)})\otimes \theta(x_{(2)})=0=\Delta(\theta(x))$, where the $x_{(1)}$ and $x_{(2)}$ live in  $\mathcal{F}_{r-1}\mathcal{H}$.
Hence, $x_{(1)}\otimes x_{(2)} \in \ker \theta \otimes \mathcal{F}_{r-1}\mathcal{H} + \mathcal{F}_{r-1}\mathcal{H} \otimes \ker \theta$, which is impossible since $r$ has been choosen to be minimal. Consequently,
$\theta$ is injective on $\mathcal{H}$, since $\mathcal{H}$ is connected.
\eproof

\begin{lemm}
\label{prim}
Let $V$ be a $K$-vector space. View $K \oplus Dipt(V)$ has a $2As^c-Dipt$-bialgebra. Then, $Prim \ (K \oplus Dipt(V))\simeq V$.
\end{lemm}
\Proof
Identify $V$ to $i(V)$ to ease notation.
Consider $(Dipt(V), \Delta)$ equipped with the semi-infinitesimal coproduct $\Delta$ of Section~4.
According to Proposition \ref{grovve}, the kernel of $\Delta$ is spanned by elements of $V$ and the:
$$\sqsubset v_1, \ldots, v_n\sqsupset_{n}= v_1 \triangleleft  (v_2 \succ (v_3 \succ  \ldots \succ(v_{(n-2)}\succ (v_{(n-1)}\succ v_n)\ldots),$$
for all $n>1$ and $v_i \in V$.
Set $X:=(v_2 \succ (v_3 \succ  \ldots \succ(v_{(n-2)}\succ (v_{(n-1)}\succ v_n)\ldots)$. Fix $n>1$. Then,
$$  \bar{\blacktriangle}(\sqsubset v_1, \ldots, v_n\sqsupset_{n}):=\bar{\blacktriangle}(v_1 \triangleleft X)=X\otimes v_1 + X_{(1)}\otimes (v_1 \triangleleft X_{(2)})\not=0,$$
because $X$ and $X_{(1)}$ has different degrees and because of the first term no linear combination of the $\sqsubset v_1, \ldots, v_n\sqsupset_{n}$ will vanish under $\bar{\blacktriangle}$, hence $\ker  \bar{\blacktriangle}\cap \ker  \Delta=V$.
\eproof

\begin{theo}
Let $\mathcal{H}$ be a $2As^c-Dipt$-bialgebra. The following are equivalent:
\begin{enumerate}
 \item{$\mathcal{H}$ is connected,}
\item{$\mathcal{H}$ is isomorphic to $ K.1_K \oplus Dipt(Prim \ \mathcal{H})$ as a $2As-Dipt$-bialgebra,}
\item{$\mathcal{H}$ is isomorphic to $K.1_K \oplus 2As^c(Prim \ \mathcal{H})$ as a 2-coassociative coalgebra.}
\end{enumerate}
Otherwise stated, the triple of operads $(2As, Dipt,Vect)$ is good.
\end{theo}
\Proof
We apply Theorem 2.3.7 \cite{GB}. Hypothesis $(H0)$ obviously holds. Hypothesis $(H1)$ is checked via Proposition \ref{2as-dipt}. Let $V$ be a $K$-vector space. As $2As^c(V)$ is the cofree 2-coassociative coalgebra over $V$, there exists a unique morphism of 2-coassociative coalgebra $\phi(V):  Dipt(V)\rightarrow 2As^c(V),$
extending the canonical map: $ Dipt(V) \twoheadrightarrow V$.
The map $\phi(V)$ is graduated and injective on $V$ since it is the identity map on $V$, hence
injective on $ Dipt(V)$ via  Lemmas \ref{prim} and  \ref{Quillen}. As for all $n>0$,
$\dim Dipt_n=\dim \ 2As_n=\dim \ 2As^c_n $, this map is an isomorphism of 2-coassociative coalgebras.
To take into account the unit, add obviously
$\phi(V)(1_K):=1_K$.
The last Hypothesis $(H2iso)$ of Theorem 2.3.7 \cite{GB} holds, hence the
triple of operads $(2As, Dipt,Vect)$ is good.
\eproof

\NB
In connected $2As^c-Dipt$-bialgebras, one has two antipodes. The usual one $S$ coming from the Hopf-algebra structure but also another one $S'$ coming from the semi-infinitesimal structure defined recursively as expected by the formula:
$$ S'(1)=1,$$
$$ \star(S'\otimes id )\vartriangle=1.\epsilon= \star(id\otimes S' )\vartriangle.$$

\section{Openings}
The motivation for this section is twofold.
Firstly, in \cite{Ltrip}, we showed the following:
\begin{theo}\cite{Ltrip}
\label{outil}
Let $\mathcal{A}$ be a binary, quadratic operad such that for any generating operations $\bullet_i \in \mathcal{A}(2)$, there exist relations only in $\mathcal{A}(3)$ of the form,
$$ (++) \ \ \ \sum_{i,j; \ \sigma_{i,j}\in \Sigma_3} \ \lambda_{i,j} \ \bullet_j( \bullet_i \otimes id) \sigma_{i,j} = \sum_{i,j; \ \sigma_{i,j}\in \Sigma_3} \ \lambda_{ij} \ \bullet_i  (id \otimes \bullet_j )\sigma_{i,j} ,$$
for any $i,j \in \{1, \ldots, \dim \mathcal{A}(2) \}$ and $\lambda_{ij} \in K$. Then, the triple $(As, \mathcal{A}, Prim_{As} \ \mathcal{A})$ endowed with nonunital infinitesimal relations is good. Considering only binary quadratic operad coming from a set operad, quadratic relations of the form:
$$(+++) \ \ \ \bullet_j( \bullet_i \otimes id)\sigma_{i,j}=(id \otimes \bullet_j )\sigma_{i,j},$$
are the only ones giving such good triples.
\end{theo}

\noindent
It is tempting to propose such a general theorem for triple of operads
endowed with the nonunital semi-infinitesimal compatibility relations instead of nonunital infinitesimal ones.

Secondly, in \cite{dipt}, the concept of L-molecule was introduced in the thesis of the author and was related to the concept coassociative covering of directed graphs. A L-molecule is a binary regular quadratic operad made out with left dipterous/ right dipterous operads whose operations $\succ_i$ and $\prec_j$
are entangled one another. Recall that two operations $\succ$ and $\prec$ are entangled if the relation:
$$ (Entanglement \ relation:) \ \ \ (x \succ y) \prec z = x \succ (y \prec z),$$
holds. For instance, the so-called predendriform operad \cite{dipt} is
a dipterous
operad $Dipt_1$ entangled with a right dipterous one $RDipt_2$ where  associative operations are merged. Here are the axioms:
$$ Dipt: \ \ (x \star_1 y)\succ_1 z = x \succ_1 (y\succ_1 z),$$
$$ (Entanglement \ relation:) \ \ (x \succ_2 y) \prec_1 z = x \succ_2 (y \prec_1 z),$$
$$ RDipt: \ \ (x \prec_2 y) \prec_2 z=  x \prec_2 (y \star_2 z), $$
and $\star_1=\star_2$.
The name molecule is borrowed from chemistry where molecules are atoms agglomerated together, the role of atoms here being played by copy of right dipterous and left dipterous operads.

These two results will be the main motivations for introducing the concept of associative molecules. First of all, we introduce the notion of dipterous like operads.

\subsection{Dipterous like operads and good triples}
For any binary regular operad $\mathcal{A}$, we set $\mathcal{A}_2:=K S$, where $S$ is the set of generating binary operations. Fix an integer $n>0$.
A dipterous like operad $^nDipt$ is a binary regular quadratic operad having a unique associative operation $\star$,
$n:=\card S$ (dipterous) operations $\succ_i$, $1\leq i \leq n$ verifying the $n$ following quadratic relations:
$$ (x \star y) \succ_i z = x \succ_i (y \succ_i z),$$
for all $i=1, \ldots,n$.
Pictorially, we represent such an operad by a circle with the associative operation inside and dipterous operations by handles:
\begin{center}
\includegraphics*[width=2cm]{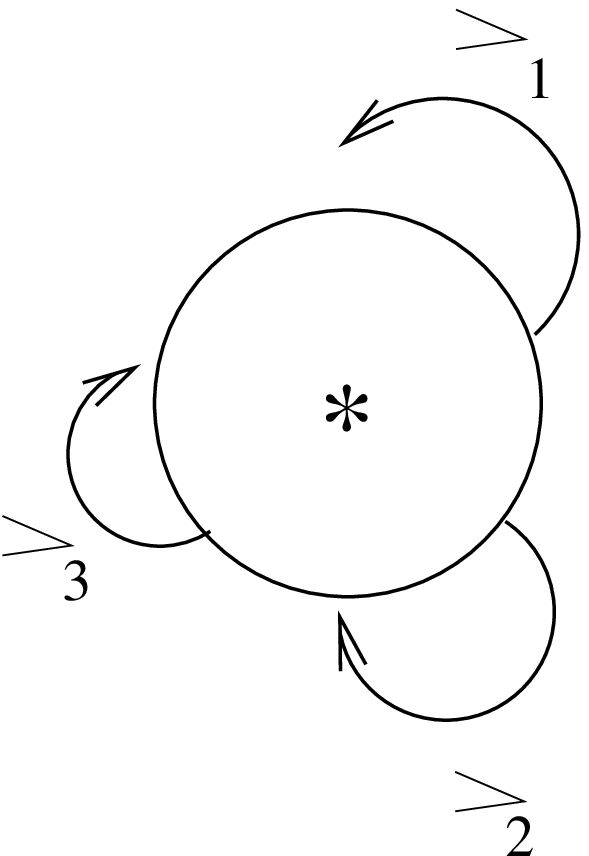}

\textsf{A dipterous like operad with 3 dipterous operations: $\succ_1, \succ_2, \succ_3$.}
\end{center}

\subsubsection{On triples of operads $(As, \ ^nDipt, Prim_{As} \ ^nDipt)$}
Infinitesimal  $^nDipt$-bialgebras are then as expected $^nDipt$-algebras equipped with a coassociative coproduct linked to operations via nonunital semi-infinitesimal relations.
\begin{theo}
Fix an integer $n>0$. The triple of operads $(As, \ ^nDipt, Prim_{As} \ ^nDipt)$ endowed with the nonunital semi-infinitesimal relations is good.
\end{theo}
\Proof
Let $V$ be a $K$-vector space.
As $^nDipt$ is a regular operad, the free $^nDipt$-algebra over $V$ will be of the form:
$$ ^nDipt(V):= \bigoplus_{m>0} \ \ ^nDipt_m \otimes V^{\otimes m},$$
the embedding map $V \hookrightarrow \ ^nDipt(V)_1 \otimes V \hookrightarrow \ ^nDipt(V)$ being denoting by $i$ (Recall that $\dim \ ^nDipt(V)_1=1$, thus $^nDipt(V)_1 \otimes V \simeq V$).
Define on $^nDipt(V)$,
the cooperation $\Delta: \ \ ^nDipt(V) \rightarrow \ \ ^nDipt(V) \otimes \ ^nDipt(V)$ recursively as follows:
$$ \Delta(i(v)):=0, \ \textrm{for all} \ v \in V,$$
$$ \Delta(x \succ_i y):= x_{(1)}\otimes (x_{(2)} \succ_i y) + (x \star y_{(1)}) \otimes y_{(2)} + x\otimes y.$$
$$ \Delta(x \star y):= x_{(1)}\otimes (x_{(2)} \star y) + (x \star y_{(1)}) \otimes y_{(2)} + x\otimes y,$$
for all $x,y \in \ ^nDipt(V)$ and $i=1, \ldots, n$.
Reapplying the same proof as in Section 2, leads to claim that $(^nDipt(V), \Delta)$
is a connected infinitesimal  $^nDipt$-bialgebra.  As
$As^c(V)$ is the cofree coassociative coalgebra over $V$, the coalgebra morphism $\phi(V): \  ^nDipt(V) \rightarrow As(V)$ extending the canonical projection $\pi': \  ^nDipt(V) \twoheadrightarrow V$ will map
any $m$-ary operations of $^nDipt_m$ to $1_K$ and will remain onto. The map
$s(V): As(V) \rightarrow \ ^nDipt(V)$, $v_1 \ldots v_m \mapsto i(v_1)\star \ldots \star  i(v_m)$ is still a coalgebra morphism and verifies $\phi(V)\circ s(V)=id$.
Theorem 2.5.1 \cite{GB} apply hence the result.

\NB
Observe, like in Section 7, that a unital version of these results can be obtained by setting $x \succ_i 1=0$ and
$1 \succ_i x=x$ for all $i=1, \ldots,n$. Use then reduced coproducts of this subsection to conclude.

\subsubsection{On triples of operads $(Com, \ ^nDipt, \ Prim_{Com }\ ^nDipt)$}
Motivated by the work of Oudom and Guin \cite{OudGuin}, another triple of operads can be interested to work with, the triple $(Com, Dipt, Prim_{Com }\ Dipt)$ endowed with the Hopf compatibility relations, where the primitive operad $Prim_{Com } \ Dipt$ has to be found. We will work with the dipterous like case. Here connected has the usual sense. To recover the definition given in \cite{GB} use the again the reduced coproduct: $\bar{\Delta}:=\Delta -( 1_K \otimes id + id \otimes 1_K).$
\begin{theo}
Fix an integer $n>0$.
If $K$ is a caracteristic zero field then
the triple of operads,
$$(Com, \ ^nDipt, \ Prim_{Com } \ ^nDipt),$$ endowed with the Hopf compatibility relations is good.
\end{theo}
\Proof
Let $V$ be a $K$-vector space. Then, the dipterous operations of $^nDipt(V)$ can be extended to $^nDipt(V)_+:=K \oplus \ ^nDipt(V)$
by requiring \cite{LodRon}:
$$ \forall x \in \ ^nDipt(V), \ \ x \star 1_K =x= 1_K \star x, $$
$$ \forall x \in \  ^nDipt(V), \ \ x \succ_i 1_K =0; \ \ \ \  1_K \succ_i x=x, $$
for all $i=1,\ldots,n$.
The symbol $1_K \succ_i 1_K$ is not defined. If $A$ and $B$ are $^nDipt$-algebras (possibly with respective units $1_A$ and $1_B$), the usual structure:
$$ (a\otimes b) \star (a' \otimes b') = (a \star_A a') \otimes (b \star_B b'),$$
$$ (a\otimes b) \succ (a' \otimes b') = (a \succ_A a') \otimes (b \succ_B b'),$$
for all $(a,b) \in A \times B$ turns the $K$-vector space
$A \otimes B$ into a  $^nDipt$-algebra
(with unit $1_A \otimes 1_B$). Hence, $^nDipt(V)_+ \otimes \  ^nDipt(V)_+$ is a $^nDipt$-algebra and the map $\delta: V \rightarrow \ ^nDipt(V)_+ \otimes \ ^nDipt(V)_+$, defined by $\delta(v)=v \otimes 1_K + 1_K \otimes v$ can be extended to a unique $^nDipt$-algebra morphism $\Delta: \ ^nDipt(V) \rightarrow \  ^nDipt(V)_+ \otimes \ ^nDipt(V)_+$ which is coassociative. This coproduct is extended to $K$ by setting $\Delta(1_K)= 1_K \otimes 1_K$. The usual flip map $\tau: \ ^nDipt(V)_+ \otimes \  ^nDipt(V)_+ \rightarrow \ ^nDipt(V)_+ \otimes ^nDipt(V)_+$, $x \otimes y \mapsto y \otimes x$ is also a morphism of dipterous like  algebras. As $\tau\delta = \delta$ so will hold $\tau\Delta= \Delta$ by unicity of the extension of the map $\delta$. Hence, $(^nDipt(V)_+, \Delta)$ is a so-called  cocommutative $^nDipt$-bialgebra, that is a $^nDipt$-algebra together with a coassociative cocommutative coproduct whose compatibility relations between operations and the coproduct is of Hopf types. It is obviously connected by construction. As $Com^c(V)$ is the cofree cocommutative coalgebra over $V$, there exists a unique morphism of coalgebras,
$$ \phi(V): \ ^nDipt(V)_+ \rightarrow Com^c(V),$$
extended the canonical projection $\pi': \
^nDipt(V)_+ \twoheadrightarrow V$ and thus verifying $\pi \circ \phi(V) = \pi'$ where $\pi: Com^c(V) \twoheadrightarrow V$ is the canonical projection onto $V$.
We have:
$$ \phi(V)(t \otimes v_1 \ldots v_m)=v_1 \ldots v_m,$$
where $t$ is any $m$-ary operation of $^nDipt_m$ and the left hand side notation stands for the symmetric tensor made over $v_1, \ldots, v_m$. Hence, $\phi(V)$ is onto.
As $\star$ is associative, the map $s(V): Com^c(V) \rightarrow \ ^nDipt(V)_+$,
$1_K \mapsto 1_K$ and $v_1 \ldots v_m \mapsto \frac{1}{m!} \sum_{\sigma \in \Sigma_m} \  i(v_1)\star \ldots \star i(v_m),$ is a coalgebra morphism such that $\phi(V)\circ s(V) =id$. Hence, Hypotheses $(H0), (H1), (H2epi)$ of Theorem 2.5.1 \cite{GB} are verified and the triple of operads $(Com, \ ^nDipt, Prim_{Com } \ ^nDipt)$ endowed with the Hopf compatibility relations is good.
\eproof
%%%%%%%%%%%%%%%%%%%%%%%%%%%%%
\subsection{On associative molecules and good triples of operads}
An associative molecule $\mathcal{A}$ is a regular quadratic binary operad made out of (maybe several copies) dipterous like operads (atoms) whose associative products are entangled. Therefore, there exist say $p$ associative products
$\star_i$ verifying:
$$ (Entanglement \ relations:) \ \  (x \star_i y)\star_j z= x  \star_i( y\star_j z),$$
for $i,j \in I$, with $I$ is any nonempty subset of $\{1,\ldots, p\}$ and for each $i\in \{1,\ldots, p\}$, there are dipterous operations $\succ_{ik}$ verifying:
$$ (x \star_i y) \succ_{ik} z = x \succ_{ik} (y \succ_{ik} z).$$
Pictorially the entanglement relation just above is represented by a dash line $i --> j$.
\begin{center}
\includegraphics*[width=7cm]{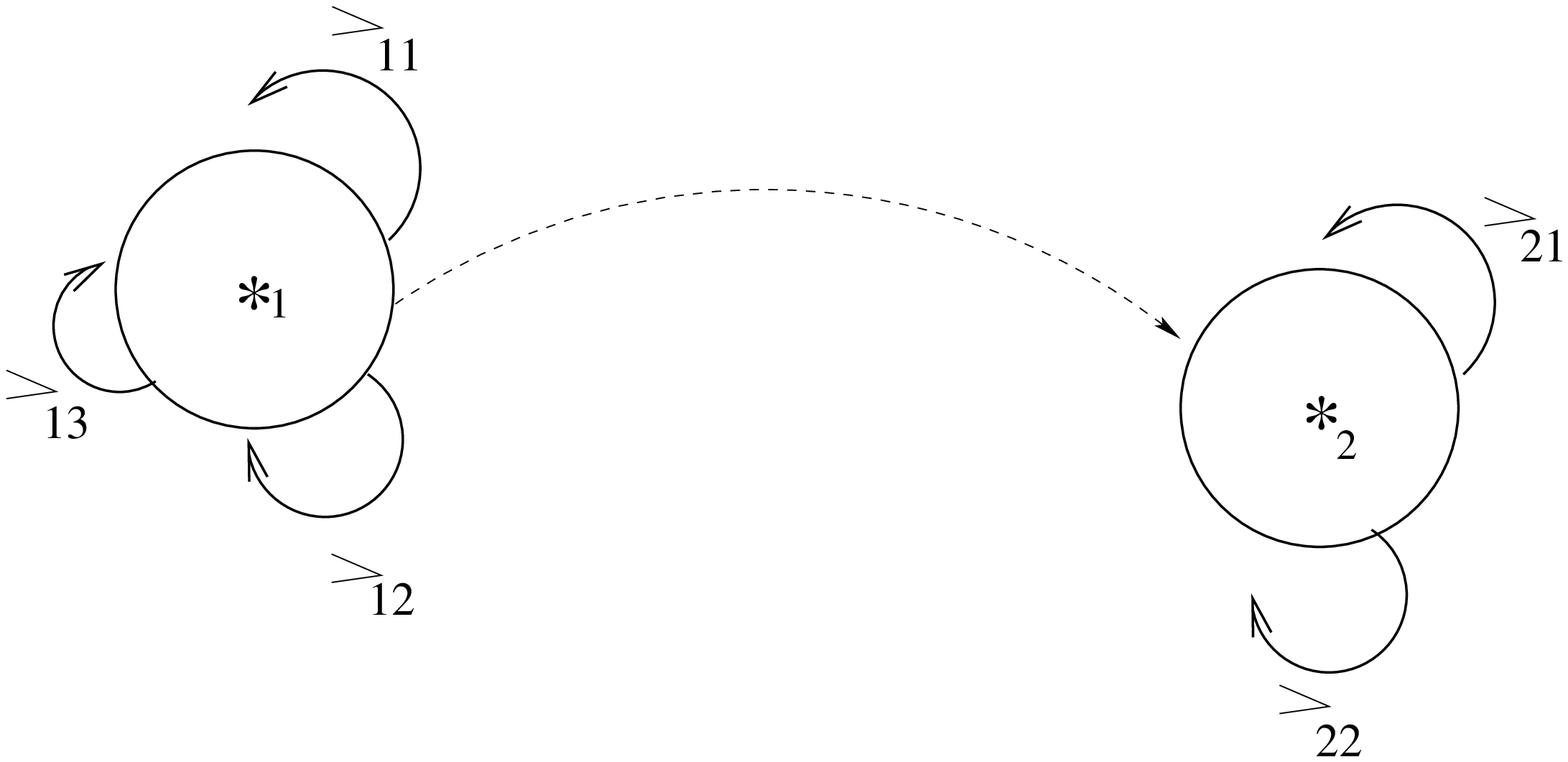}

\textsf{An associative molecule: Entanglement of two dipterous like operads.}
\end{center}
\noindent
As nonunital semi-infinitesimal relations coincide with nonunital  infinitesimal for associative products,
the same arguments of Subsubsection 7.1.1 apply and any associative molecules $\mathcal{A}$ will lead to good triples of operads $(As, \mathcal{A}, Prim_{As} \ \mathcal{A})$ endowed with the semi-infinitesimal relations.
\NB
The same holds by replacing the adjective dipterous by right dipterous.
\NB
In this subsection, one is forced to use only the nonunital semi-infinitesimal relations since the action of the unit does not behave well with the entanglement relations $ (x \star_i y)\star_j z= x  \star_i( y\star_j z)$.

\bibliographystyle{plain}
\bibliography{These}

\end{document}